\documentclass{amsart}
\usepackage{amssymb,latexsym}
\usepackage[dvips]{graphicx}
\newtheorem{theorem}{Theorem}[section]
\newtheorem{prop}[theorem]{Proposition}
\newtheorem{lemma}[theorem]{Lemma}
\newtheorem{remark}[theorem]{Remark}
\newtheorem{question}[theorem]{Question}
\newtheorem{definition}[theorem]{Definition}
\newtheorem{cor}[theorem]{Corollary}

\DeclareMathOperator{\ep}{\epsilon}
\def\co{\colon\thinspace}
\begin{document}

\title{Symplectic forms and surfaces
of negative square}

\author{Tian-Jun Li \& Michael Usher}
\address{School  of Mathematics\\  University of Minnesota\\ Minneapolis, MN 55455}
\email{tjli@math.umn.edu}
\address{Department of Mathematics\\ Princeton University\\ Princeton, NJ  08540}\email{musher@math.princeton.edu}

\begin{abstract}
We introduce an analogue of the inflation technique of
Lalonde-McDuff, allowing us to obtain new symplectic forms from
symplectic surfaces of negative self-intersection in symplectic
four-manifolds.  We consider the implications of this construction
for the symplectic cones of K\"ahler surfaces, proving along the way
a result which can be used to simplify the intersections of distinct
pseudoholomorphic curves via a perturbation.
\end{abstract}

\thanks{The first author is supported in part by NSF grant 0435099 and the McKnight fellowship.
The second author is supported by an NSF postdoctoral fellowship.}
\maketitle

\section{ Introduction}

Given an embedded symplectic surface $C$ of non-negative
self-intersection in a symplectic 4-manifold $(M,\omega)$, the
inflation process in \cite{LM} gives rise to new symplectic forms in
the  class $[\omega]+tPD[C]$ for arbitrary $t>0$. In this paper we
show that there is an analogous construction in the case of an
embedded symplectic surface of negative self-intersection.

\begin{theorem} \label{main}
Suppose $C$ is an embedded connected symplectic surface representing
a class $e$ with $e\cdot e=-k<0$ and $a=\omega(e)$. Let $h=k$  if
$C$ has positive genus or $C$ is a sphere with $k$ even, and $h=k+1$
if $C$ is a sphere with $k$ odd. Then there are symplectic forms
$\omega_t$ representing the classes $[\omega] +tPD(e)$ for any $t\in
[0, {2a\over h})$.
\end{theorem}

This is achieved by the normal connected sum construction (see
\cite{Go}, \cite{MW}). In fact the inflation process can be viewed
this way as well. However there are two distinct features from the
inflation process. The first is the upper bound on $t$. The second
is that the surface $C$ 
is not symplectic with respect to
the forms $\omega_t$ when $t > {a\over k}$ 
(such values of $t$ occur as long as $C$ is not a ($-1$)-sphere).
Indeed, the symplectic area of $C$ is non-positive for these values
of $t$.

From the known characterization of the symplectic cones of $S^2-$bundles
\cite{Mc3}, for any triple $(g,k,a)$ with $g\geq 0$ and
$k,a>0$ and any $\ep>0$ it is a routine exercise to find a
symplectic $4$-manifold $(M,\omega)$ containing a symplectic surface
$\Sigma_{\ep}$ of genus $g$, square $-k$, and area $a$ such that
$[\omega]+(2a/h+\ep)PD[\Sigma_{\ep}]$ does \emph{not} admit
symplectic forms, where $h$ is as in the statement of Theorem
\ref{main}. In this regard, Theorem \ref{main} may be considered a
best possible result for the generality in which it is stated.

Using the pairwise normal connected sum construction
we 
will also show how to apply the 
construction of Theorem \ref{main} to a configuration of surfaces
intersecting each other positively and transversally.

To  apply such a construction we need  to locate configurations of
surfaces. Such configurations sometimes appear as pseudo-holomorphic
curves. It is shown in \cite{Mc4} that any irreducible simple
pseudo-holomorphic curve can be perturbed to a pseudo-holomorphic
immersion, possibly after a $C^1$-small change in the almost complex
structure. We show how to further perturb such an immersion to an
embedding. In fact we are able to show that any configuration of
simple $J$-holomorphic curves can be perturbed to a configuration of
symplectic surfaces which intersect each other positively and
transversally and which are $J'$-holomorphic for an almost complex
structure arbitrarily $C^1$-close to $J$.

Holomorphic curves of negative self-intersection actually
characterize the K\"ahler cone by (the extension of) the
Nakai-Moishezon criterion. Thus it is interesting to apply this
construction to  K\"ahler surfaces. Let $(M,J)$ be a K\"ahler
surface and $H_J^{1,1}$ denote the real part of the $(1,1)-$subspace
of $H^2(M;{\mathbb C})$ determined by $J$.  The classical Hodge
index theorem then asserts that the restriction of the intersection
form to $H^{1,1}_{J}$ is a bilinear form of type $(1,h^{1,1}-1)$.
The \emph{positive cone} of $H^{1,1}_{J}$ is then by definition the
set of classes in $H^{1,1}_{J}$ which have positive square and pair
positively with the class of the given K\"ahler form. Buchdahl and
Lamari have recently independently proven the following result:
\begin{theorem}\label{bl}\emph{(}\cite{B},\cite{L}\emph{)} For a K\"ahler
surface $(M,J)$, any class in the positive cone of $H_{J}^{1,1}$  is
represented by a K\"ahler form if it is positive on each holomorphic
curve with negative self-intersection. \end{theorem} (Note that the
Hodge index theorem implies that any class in the positive cone of
$H^{1,1}_{J}$ is automatically positive on each curve of
\emph{non}-negative self-intersection.) Applying Theorem \ref{main}
to a curve $C$ of negative self-intersection, the K\"ahler cone can
be enlarged across the ``wall'' consisting of cohomology classes
which vanish on $[C]$ unless $C$ is a $(-1)$-sphere. This suggests
the following symplectic Nakai-Moishezon criterion:

\begin{question}\label{q} For a  K\"ahler surface $(M,J)$, is
every class in the positive cone of  $H^{1,1}_{J}$ which is positive
on each $(-1)$-sphere (possibly reducible) represented by a symplectic form?
\end{question}

To motivate this, note that by the Riemann-Roch theorem and the
adjunction formula the expected dimension of the space of embedded
pseudoholomorphic genus $g$ curves in the class $[C]$ is \[
d([C])=2(g-1+\langle c_1(M),[C]\rangle)=[C]\cap [C]-\langle
c_1(M),[C]\rangle=[C]\cap [C]+1-g,\] which as the last expression
above demonstrates is negative if $[C]$ is the class of any
negatively self-intersecting curve other than a $(-1)$-sphere. Thus
for generic almost complex structures $J'$ close to $J$, there will
be no $J'$-holomorphic curves in the class $C$. The theory of
pseudoholomorphic curves hence does not provide any obstruction to
deforming the symplectic form to one which pairs negatively with
$C$. If $C$ is the class of a $(-1)$-sphere, on the other hand,
Gromov-Taubes theory shows that any symplectic form deformation
equivalent to the K\"ahler form must pair positively with $C$.

 When $p_g=0$  Question \ref{q}  has an affirmative answer.
In this case  we have $b^+=1$, so every class of positive square
which is positive on $-1$ symplectic spheres is realized by a
symplectic form (\cite{LL}). In addition, for a minimal surface of
general type, the canonical class $K$ has been shown to be in the
symplectic cone (\cite{S}, \cite{Ca}).

In a more general setting, the answer to Question \ref{q} seems
elusive.  Our methods do, however, enable us to progress somewhat
farther on the following related question:
\begin{question}\label{q2} Let $\{C_1,\ldots,C_n\}$ be reduced
irreducible holomorphic curves of negative square, none of which is
a $(-1)$-sphere, such that there exist classes $\alpha$ in the
positive cone of $H^{1,1}_{J}$ satisfying $\langle \alpha,
[C_i]\rangle<0$ for each $i$.  Do some of these classes $\alpha$
admit symplectic forms?
\end{question}

In Section 4 we outline methods for using Theorem \ref{main} to
answer this question in certain situations, and we illustrate these
methods by applying them in detail to all of the subsets of a
particular set of 21 negative-square curves in a rigid surface $K$
that was introduced in \cite{KK}.  We choose a rigid surface as our
primary example in order to ensure that the curves in question
cannot be made to disappear by an integrable variation in the
complex structure; as such, we may state with certainty that the new
symplectic forms that we construct are not directly obtainable by
considerations of K\"ahler geometry.

The methods of Section 4 can be applied to a wide variety of
configurations of the curves $C_1,\ldots,C_n$ in Question \ref{q2},
but there are also many configurations to which these methods do not
apply.  It seems unlikely that there is any necessary and sufficient
condition on the configuration that can be expressed at all
concisely, but we provide an example of a moderately general
sufficient condition in Theorem \ref{ADE}.



We would like to thank D. McDuff for her valuable suggestions on
how to extend her result in \cite{Mc4} to our situation. The first author is also grateful to Y. Ruan for discussions on the
6-dimensional symplectic minimal model program which inspired Theorem \ref{main}.

\section{The construction}

Theorem \ref{main} is an application of  the normal connected sum
construction with symplectic $S^2-$bundles. So let us collect some
facts about symplectic structures on such manifolds and embedded
symplectic surfaces in them.

Up to diffeomorphisms, there are two orientable $S^2-$bundles over a Riemann surface $\Sigma$: the trivial one $\Sigma\times S^2$,
and the non-trivial one $M_{\Sigma}$. By \cite{LM}, symplectic forms on $S^2-$bundles are determined by their cohomology classes up
to isotopy. Thus we can pick any convenient symplectic form in a fixed cohomology class.

We begin with the easier case: the product bundle. In this case we
use split forms as our model forms. Clearly every class of the
positive cone is represented by a split symplectic form. And for a
split symplectic form, the vertical fibers and  horizontal sections
are symplectic. The class of any  section with (even) positive
square is then represented by an immersed symplectic surface with
only positive transverse self-intersections, which can then be
smoothed to an embedded symplectic surface.

Now let us deal with the non-trivial bundle $M_{\Sigma}$ over a
positive genus surface. We use K\"ahler forms as our model forms.
The following result is essentially contained in \cite{Mc3} and
\cite{H} (we present it here since it may not be very well-known).

\begin{prop} \label{stablebundle} Let ${\mathcal E}$ be a holomorphic rank 2 bundle over $\Sigma$ with  $g(\Sigma)>0$ and $c_1({\mathcal E})=-1$. Let
$(M,J_{\mathcal E})$ be the complex ruled surface $P({\mathcal E})$.
Then the K\"ahler cone is the positive cone if and only if ${\mathcal
E}$ is stable. Furthermore, for appropriately chosen holomorphic
structures on ${\mathcal E}$, the class of any section with (odd)
positive square can be represented by an embedded surface which is
symplectic with respect to any K\"ahler form.
\end{prop}

\begin{proof} Notice that the slope of ${\mathcal E}$ is $-{1\over 2}$. Therefore the stability of ${\mathcal E}$
is equivalent to the statement that every holomorphic line subbundle
${\mathcal L}$ of ${\mathcal E}$ has $c_1({\mathcal L})\leq -1$.
Observe that any holomorphic line bundle ${\mathcal L}\subset
{\mathcal E}$ gives rise a to a holomorphic section $Z({\mathcal
L})$ of $P({\mathcal E})$, and vice versa. Since the normal bundle
to $Z({\mathcal L})$ is ${\mathcal L}^*\otimes {\mathcal
E}/{\mathcal L}$, all sections of $P({\mathcal E})$ have positive
self-intersection if and only if $\mathcal{E}$ is stable.
 The statement about K\"ahler cone now follows from the arguments in Proposition 3.1 in \cite{Mc3}
 (see also \cite{H}).

For the second statement,  it suffices to show that the class
$[s^+]$ of a section with square $+1$ is symplectic. As all the
fibers are holomorphic and hence symplectic and the classes of
sections with higher squares have form $[s^+]+m[fiber]$ for $m>0$,
these classes are  represented by positively immersed symplectic
surfaces, which can be smoothed to embedded ones. We may take the
holomorphic structure on ${\mathcal E}$ to be that on a non-trivial
extension of ${\mathcal L}$ by the trivial line bundle ${\mathcal
O}$, where ${\mathcal L}$ is a degree $-1$ holomorphic line bundle.
The section $Z({\mathcal L})$ is then a holomorphic, and so in
particular symplectic, $+1$ section.
\end{proof}

Finally, the non-trivial bundle over a sphere is diffeomorphic to
the blow up of ${\mathbb CP}^2$ at a point, and the exceptional
divisor is a section with square $-1$. As is well-known, either
using the standard symplectic reduction picture or  algebraic
geometry, we can construct symplectic forms in every class in the
positive cone which is positive on the class of a section with square
$-1$, such that, for every odd $k\geq -1$,  there are symplectic
sections with square $k$.

Now we are ready to  prove Theorem \ref{main}.
\begin{proof} Let $R$ be the  trivial sphere bundle over
the surface of genus $g(C)$ if $k$ is even, and the non-trivial one
if $k$ is odd. Let $s^{\pm k}$ be the class of a section with square
$\pm k$. Then $s^{+k}$ and $s^{-k}$ form a basis for $H_2(R;{\mathbb
Z})$. Since $s^{+k}\cdot s^{-k}=0$, a cohomology class of the form
$$c^+PD(s^{+k})+c^-PD(s^{-k})$$ has positive square if and only if $c^+>|c^-|$.

Suppose first that $C$ is not a sphere with $k$ odd. By Proposition
\ref{stablebundle} and the discussions preceding it,
 there exists a
 symplectic form $\tau_t$ on $R$ in
the class
$${a\over
k}PD(s^{+k})+(t-{a\over k})PD(s^{-k})$$
 for any $t\in (0, {2a\over
k})$.
 By Proposition
\ref{stablebundle}, there is a $\tau_t-$symplectic section $S^{+k}$
in the class $s^{+k}$. Notice that the symplectic surfaces $C$ and
$S^{+k}$ have opposite self-intersection and equal symplectic area
$a$.
 Thus we can perform the symplectic sum construction
to $(M,C,\omega)$ and  $(R, S^{+k},\tau_t)$ to obtain a new
symplectic manifold $(X,\omega_t)$. As observed in \cite{Go},  $X$
and $M$ are diffeomorphic. Moreover,
 because the surface
$S^{-k}$ is disjoint from surface $S^{+k}$ in $R$,  it becomes a
surface in $M$ which is homologous to $C$.  Thus we have
$$\omega_t(e)=\tau_t(s^{-k})=a-tk.$$
Therefore $[\omega_t]=[\omega]+tPD(e)$.

In the case that $C$ is  a sphere with $k$ odd, by the discussions
after Proposition \ref{stablebundle}, there exists a symplectic form
in the class $c^+ PD(s^{+1})+ c^- PD(s^{-1})$ if and only if
$c^+>-c^->0$. We would like to express the condition  in terms of
the basis $s^{+k}$ and $s^{-k}$. Since $s^m\cdot s^n={m+n\over 2}$,
we see that a class $\alpha$ contains a symplectic form if and only
if
$${\alpha(s^{-k})\over \alpha(s^{+k})}>{1-k\over k+1}.$$
 Thus the allowed values of $t$ are those between $0$ and ${2a\over
k+1}$, as claimed.
\end{proof}

\begin{remark} Notice that
$$\omega_t\cdot \omega_t=\omega^2+2t\omega(e)+t^2e\cdot e=-t^2k+2ta +\omega^2=\omega^2+k[{a^2\over k^2}-(t-{a\over k})^2].$$
So the volume of the symplectic manifold $(M,\omega_t)$ is greater
than that of $(M,\omega)$ for each $t\in (0, 2a/k)$.
\end{remark}

We can generalize Theorem \ref{main} to a
configuration of transversally intersecting symplectic surfaces.

\begin{theorem} \label{configuration} Suppose $C_1,..., C_l$ is a set of embedded
symplectic surfaces with self-intersection $C_i\cdot C_i=-k_i<0$ and
intersecting positively and transversally.  Let $e_i$ be the class
of $C_i$. Then there are symplectic forms in the class
$[\omega]+\sum_i t_i PD(e_i)$ for any $0<t_i<2\omega(e_i)/h_i$,
where $h_i$ is as in Theorem \ref{main}.
\end{theorem}

\begin{proof} We prove the theorem in the case that there are only two curves $C_1$ and $C_2$.
The idea for the general case is the same.

 The key
point is that the symplectic sum construction in Theorem \ref{main},
when applied to $C_1$,  can be done in a way such that $C_2$,
possibly after an isotopy, is still symplectic with respect to the
new symplectic structures $\omega_1$, which is in the class
$[\omega]+t_1PD(e_1)$ with $t_1\in (0, 2\omega(e_i)/h_1)$. This is
possible due to the pairwise sum feature in \cite{Go}.

First, by applying Lemma 2.3 of \cite{Go},  perturb $C_2$  such that
$C_2$ intersects $C_1$ orthogonally with respect to $\omega$. Since
the fiber spheres in $R$ are symplectic and intersect the symplectic
section $S^{+k_1}$ transversally, we can likewise assume that the
symplectic section $S^{+k_1}$ intersects a total of $k_{12}=C_1\cdot
C_2$ fibers, all orthogonally. Denote this union of the fibers by
$F$. Now apply pairwise sum to $(M, C_1, C_2)$ to $(R, S^{+k_1}, F)$
to get a symplectic surface $C_2'$.  Finally, apply the symplectic
sum construction to $C_2'$ as in the proof of Theorem \ref{main}.

\end{proof}

\begin{remark} Notice that since $e_1\cdot e_i\geq 0$ for
$i\geq 2$,  $[\omega]+t_iPD(e_1)$ is positive on $e_2$ for $t_1$
positive. One has $S_i^+=S_i^-+k_if$ where $f$ is the homology class
of the fiber in $R$, so
$$\int_{S_i^-}\tau=\int_{C_i}\omega+t_iPD(e_i)=a_i-t_ik_i=\int_{S_i^+}\tau-k_i\tau(f).$$
Thus $\tau(f)=t_i$. This is consistent with the normal connected sum
picture. The area of the surface $C_j$ increases by $(e_j\cdot
e_i)\tau(f)$, which is indeed equal to $t_i(e_i\cdot e_j)$.
\end{remark}

\begin{remark} If these surfaces actually intersect, then some of the values of $t_i$ can
be taken larger than in the statement of the theorem.
\end{remark}

\section{Configurations of embedded symplectic surfaces and  pseudo-holomorphic curves}

In attempting to answer questions such as Question \ref{q}, we might
wish to apply Theorem \ref{configuration} to some finite set of
holomorphic curves.  However, the proof of Theorem
\ref{configuration} depends on the assumption that the symplectic
submanifolds being considered intersect positively and transversely,
which is a property that our set of holomorphic curves might not be
known to have.  Assume that we are given a collection of distinct
$J$-holomorphic curves $C_1,\ldots,C_k$ in the symplectic
$4$-manifold $M$ (we adopt the convention that a $J$-holomorphic
curve is the image of a generically injective $J$-holomorphic map
from some irreducible compact Riemann surface). Corollary 4.2.1 of
\cite{Mc4} asserts that, at the possible cost of $C^1$-slightly
changing the almost complex structure $J$, we may perturb any one of
these curves to a pseudo-holomorphic immersion.  We first give a
simple modification of McDuff's argument to show that, in fact, we
may perturb all of the curves and the almost complex structure
simultaneously so that $C_1,\ldots,C_k$ become immersed.

\begin{lemma} \label{immerse}  Let $u_i\co \Sigma_i\to M$ be $J$-holomorphic maps
with images $C_i$.  Then given $\ep >0$ there are an almost complex
structure $\tilde{J}$ and $\tilde{J}$-holomorphic immersions
$\tilde{u}_i\co\Sigma_i\to M$ such that
$\|\tilde{u}_i-u_i\|_{C^2}<\ep$ and
$\|\tilde{J}-J\|_{C^1}<\epsilon$.\end{lemma}
\begin{proof}
Let $p\in M$ be a critical value for one or more of the $u_i$.  It
is shown in \cite{Mc4} that the various $u_i$ each have just
finitely many critical points, so denote the various critical points
in $\cup \Sigma_i$ having image $p$ by $z_1,\ldots, z_m$.  For
$j=1,\ldots,m$, if $z_j\in \Sigma_i$ let $D_j\subset \Sigma_i$ be a
disc around $z_j$, and let $v_j=u_i|_{D_j}$.  By shrinking the
various $D_j$, we assume that the $D_j$ are disjoint and that $z_j$
is the only critical point of the restriction $v_j$.  Since the
intersections (and self-intersections) of the various $C_i$ are
isolated, let $U\subset M$ be a coordinate neighborhood of $p$ in
which the $C_i$ meet each other and themselves only at $p$ and such
that for each $j$ $v_{j}^{-1}(U)\subset D_j$ is a connected
component of $\cup u_{i}^{-1}(U)$.  Shrinking $U$ if necessary,
assume also that $U$ contains no critical values of the various
$u_i$ other than $p$. Now fix neighborhoods $W_m\subset
U_m\subset\cdots\subset W_1\subset U_1\subset U$ of $p$.   By
Theorem 4.1.1 of \cite{Mc4}, there is a family $v_{1}^{\delta}$
($\delta>0$) of $J$-holomorphic immersions $D_1\to M$, converging in
$C^2$ norm to $v_1$ as $\delta\to 0$.  For $\delta$ small, define
$\tilde{v}_1(z)=\chi(z) v_{1}^{\delta}(z)+(1-\chi)(z)v_1(z)$, where
$\chi$ is a smooth cutoff function which is $1$ on a neighborhood of
$v_{1}^{-1}(W_1)$ and $0$ on a neighborhood of the complement of
$v_{1}^{-1}(U_1)$. $\tilde{v}_1$ is then $C^2$-close to $v_1$, so
there is an almost complex structure $J'_1$ which agrees with $J$
away from $U_1\setminus W_1$, makes $\tilde{v}_1$ $J'_1$
holomorphic,
 and is $C^1$-close to $J$ everywhere (see the proof of Corollary 4.2.1 of \cite{Mc4}, or the proof of Proposition \ref{pert} below).
Furthermore,  if $U\cap Im v_1$ is a distance at least $K$ from
$\cup C_i\setminus Im v_1$, then for $\delta$ small enough $U\cap Im
\tilde{v}_1$ will be a distance $K/2$ from $\cup C_i\setminus Im
v_1$, and so using a cutoff function supported in a
$(K/3)$-neighborhood of $Im (\tilde{v}_1)\cap (U_1\setminus W_1)$,
we can patch together $J$ and $J'_1$ to obtain an almost complex
structure $\tilde{J}_1$ which is $C^1$-close to $J$, agrees with $J$
outside $U_1\setminus W_1$ and on a neighborhood of $\cup
C_i\setminus Im(\tilde{v}_1)$, and makes $\tilde{v}_1$
pseudo-holomorphic.

With this done, we now apply the same procedure sequentially to
$v_2,\ldots,v_m$, obtaining almost complex structures $\tilde{J}_j$
which are $C^1$-close to $J$ globally and which agree with
$\tilde{J}_{j-1}$ both near $\cup C_i\setminus Im(v_j)$ and outside
$U_j\setminus W_j$, and $\tilde{J}_{j}$-holomorphic immersions
$\tilde{v}_j$ which are $C^2$-close to $v_j$.  Modifying the
original maps $u_i\co \Sigma_i\to M$ by replacing the restrictions
$v_j\co D_j\to M$ by $\tilde{v}_j$, we get $\tilde{J}_m$-holomorphic
maps $\tilde{u}_i$ which have no critical values inside $U$ and
agree with the $u_i$ outside $U$.  So we have reduced the number of
critical values by 1, and repeating the process at each critical
value gives the almost complex structure $\tilde{J}$ and the
$\tilde{J}$-holomorphic immersions $\tilde{u}_i$ that we desire.
\end{proof}

Applying Lemma \ref{immerse}, we may assume that we now have a set
of distinct immersed $J$-holomorphic curves $C_i$, and we aim now to
show that these curves may be perturbed further to a set of
symplectic surfaces $C'_i$ whose intersections are all transverse
and positive with $C'_i\cap C'_j\cap C'_k=\varnothing$ when
$i$,$j$,$k$ are all distinct.  In fact, our perturbed curves $C'_i$
will agree with $C_i$ outside an arbitrarily small neighborhood of
the initial intersection points; will be arbitrarily $C^1$-close to
$C_i$ (from which it immediately follows that they are symplectic);
and will be made simultaneously pseudoholomorphic by an almost
complex structure $J'$ arbitrarily $C^1$-close to $J$.

We start by finding a nice coordinate system near any given
intersection point of our curves.  In the case where $J$ is
integrable, any given holomorphic coordinate chart can be modified
by an element of $GL(2,\mathbb{C})$ to satisfy the conditions we
need, so the arguments below are only needed in the non-integrable
case.

\begin{lemma}\label{coord}
Given immersed $J$-holomorphic curves $C_0,\ldots,C_m\subset M$ all
having an isolated intersection at the point $p$, there is a
coordinate chart $U$ around $p$ with coordinates $z,w$ such that:
\begin{itemize}
\item[(i)] $C_0\cap U=\{(z,w)\in U|w=0\}$, \item[(ii)] Each set
$\{(z,w)\in U|w=const\}$ is $J$-holomorphic, and \item[(iii)] For
$i\geq 1$ there are smooth functions $g_i$ of the form
$g_i(z)=a_iz^{k_i}+O(|z|^{k_i+1})$ $(a_i\neq 0$,$k_i\geq 1)$ such
that $C_i\cap U=\{(z,g_i(z))\}$
\end{itemize}
\end{lemma}
\begin{proof}  A coordinate chart $U_0=\{(z',w')\}$ satisfying (i) and
(ii) may be constructed by using Lemmas 5.4 and 5.5(d) of \cite{Ta}.
To obtain condition (iii), first note that for a generic linear
change of coordinates $(z',w')\mapsto (z'+cw',w')$ we retain
properties (i) and (ii) and additionally ensure that $\{z'=0\}$ is
transverse to each of the $C_i$.  Now condition (ii) implies that
the antiholomorphic tangent space of our  almost complex structure
$J$ is given in these coordinates by
\[
T^{0,1}_{J}=\langle
\partial_{\bar{z'}}+\alpha(z',w')\partial_{z'},v(z',w')\rangle\] for a
certain function $\alpha$ and complexified vector field $v$. By
Ahlfors-Bers' Riemann mapping theorem with smooth dependence on
parameters \cite{AB}, the equation
$u_{\bar{z'}}+\alpha(z',w')u_{z'}=0$ can be solved for a smooth
function $u(z',w')$ with $u(0,w')=0$.

Changing coordinates to $(z,w)=(z'+u(z',w'),w')$, we have that
$\{z=0\}=\{z'=0\}$ is transverse to each of the $C_i$, so that
after possibly shrinking the coordinate chart $U$ we have $C_i\cap
U=\{(z,g_i(z))\}$ for some smooth functions $g_i$.  In terms of
the coordinates $(z,w)$, we have for certain functions $a$ and $b$
both vanishing at the origin,\[ T^{0,1}_{J}=\langle
\partial_{\bar{z}},\partial_{\bar{w}}+a(z,w)\partial_w+b(z,w)\partial_z\rangle.\]
It is then a simple matter to check that a curve
$\{(z,g(z))\}\subset U$ is $J$-holomorphic exactly
if\[b(z,g(z))=\frac{g_{\bar{z}}-a(z,g(z))g_z}{|g_{\bar{z}}|^2-|g_z|^2}.\]
But then the fact that $a(z,g(z))$ and $b(z,g(z))$ are smooth
functions of $z$ and vanish at $z=0$ implies that the lowest-order
terms in the Taylor expansion of $g$ are functions only of $z$ and
not of $\bar{z}$.  Of course, our functions $g_i$ can't be
constants (since the $C_i$ ($i\geq 1$) have an isolated
intersection point with $C_0=\{w=0\}$ at the origin), so it
follows that the $g_i$ all have the form specified in condition
(iii).\end{proof}

\begin{prop} \label{pert} In the situation of Lemma \ref{coord}, given any
sufficiently small $\delta>0$ there is a surface $C^{\delta}_{0}$
such that, where $B_{\delta}=\{(z,w)\in U||z|< \delta\}$,
$C^{\delta}_{0}\cap(X\setminus B_{\delta})=C_0\cap(X\setminus
B_{\delta})$, while all intersection points of $C^{\delta}_{0}$
with $C_i$ ($i>0$) that are contained in $B_{\delta}$ are in fact
contained in $B_{\delta^2}$ and are transverse, positive, and
distinct from $p$ and from each other as $i$ varies.  Further
there is a constant $A$ depending on the curves $C_i$ but not on
$\delta$ such that $dist_{C^2}(C^{\delta}_{0},C_0)\leq A\delta^2$,
and there is an almost complex structure $J'$ agreeing with $J$
near $C_i$ ($i>0$) and making $C^{\delta}_{0}$ holomorphic with
$\|J'-J\|_{C^1}\leq A\delta^2$.\end{prop}
\begin{proof}
Work in coordinates provided by the conclusion of Lemma
\ref{coord}, and let $c_i$ be constants such that for $i>0$ \[
|g_i(z)-a_iz^{k_i}|<c_i|z^{k_i+1}|.\] Given $\ep>0$ write
\[R_{\ep}=\max_{i\geq 1}\left(\frac{2\ep}{|a_i|}\right)^{1/k_i};\]
we will only work with $\ep$ so small that \[
\sqrt{R_{\ep}}<\min_{i\geq 1}\frac{|a_i|}{2c_i}\]  Then for any
such $\ep$, if $R_{\ep}\leq |z| \leq \sqrt{R_{\ep}}$, we have\[
|g_i(z)|\geq (|a_i|-c_i|z|)|z|^{k_i}>\frac{|a_i|}{2}|z|^{k_i}\geq
\ep.\]

Write $\delta=\sqrt{R_{\ep}}$; note that we may alternatively
express $\ep$ in terms of an arbitrary $\delta >0$, and then for
$\delta$ small enough $\ep$ is bounded by a constant times
$\delta^2$. Fix a cutoff function $\chi(z)$ with image $[0,1]$ equal
to one for $|z|\leq \delta^2$ and zero for $|z|\geq\delta$, with
$\|\chi\|_{C^2}<\frac{4}{\delta^2}$. Let
$C_{0}^{\delta}=\{(z,\ep^2\chi(z))\}$; obviously $C^{\delta}_{0}$
agrees with $C_0$ outside $B_{\delta}$. Since $\ep^2\chi(z)\leq
\ep^2<\ep$ while each $|g_i(z)|>\ep$ for $|z|\in[\delta^2,\delta]$,
evidently the intersection points of $C_i$ with $C^{\delta}_{0}$
contained in $B_{\delta}$ are just those points $(z,\ep^2)$ with
$|z|<\delta^2=R_{\ep}$ such that $g_{i}(z)=\ep^2$.

Write $\tilde{g}_i(z)=\ep^{-2}g_i\left(
\left(\frac{\ep^2}{a_i}\right)^{1/k_i}z\right)$; then
\[\tilde{g}_i(z)=1\Leftrightarrow
g_i\left(\left(\frac{\ep^2}{a_i}\right)^{1/k_i}z\right)=\ep^2,\] so
the intersections of $C_{i}$ with $C^{\delta}_{0}$ are of just the
same type as the intersections of the graph of $\tilde{g}_i(z)$ with
$\{w=1\}$. Now we see that $\tilde{g}_i(z)=z^{k_i}+\tilde{r}_i(z)$
where $|\tilde{r}_i(z)|\leq\tilde{c}_i\ep^{2/k_i}|z|^{k_i+1}$. Hence
the graph of $\tilde{g}_i(z)$ is $O(\ep^{2/k_i})$ away in $C^1$ norm
from that of $z\mapsto z^{k_i}$, so since the latter's only
intersections with $\{w=1\}$ are positive and transverse at the
$k_i$th roots of unity, for $\ep$ small enough the graph of
$\tilde{g_i}$ will also have just $k_i$ distinct positive transverse
intersections with $\{w=1\}$, each at a point a distance
$O(\ep^{2/k_i})$ from a different one of the $k_i$th roots of unity.
Scaling back, we conclude that the intersections of $C_i$ with
$C^{\delta}_{0}$ that are contained in $B_{\delta}$ are in fact
contained in $B_{\delta^2}$ and are transverse, positive, and
located at points a distance $O(\ep^{4/k_i})$ from the various
$(\ep^2/a_i)^{1/k_i}\eta$ for $\eta$ a $k_i$th root of unity.

Obviously for any given $i$ the points of $C_i\cap C^{\delta}_{0}$
are all distinct for small enough $\ep$.  For small enough $\ep$
these intersections  vary continuously in $\ep>0$, so if it
weren't the case that the sets $C_i\cap C_j\cap C^{\delta}_{0}$
were all eventually empty for $\ep$ small enough and $i,j$
distinct, we would then, by varying $\ep$, obtain a continuous
family of points in $C_i\cap C_j$, which is impossible since $C_i$
and $C_j$ are distinct holomorphic curves and so have isolated
intersections.

Finally, note that $\|\ep^2\chi\|_{C^2}\leq\ep^2(4/\delta^2)\leq
A\delta^2$ for a certain constant  $A$ and $\delta$ sufficiently
small, so that $C^{\delta}_{0}$ is indeed less than $A\delta^2$
away from $C_0$ in $C^2$ norm.  Letting $\beta (w)$ be a cutoff
function which is $1$ for $|w|<2\ep^2$ and $0$ for $|w|\geq \ep$,
if $J$ is defined by $T^{0,1}_{J}=\langle
\partial_{\bar{z}},v\rangle$ then setting \[T^{0,1}_{J'}=\langle
\partial_{\bar{z}}+\beta(w)\left((\ep^2\chi)_{\bar{z}}\partial_w+(\overline{\ep^2\chi})_{\bar{z}}\right)\partial_{\bar{w}},v\rangle\]
defines an almost complex structure $J'$ which makes
$C^{\delta}_{0}$ holomorphic and which (since $|g_i(z)|>\ep$
whenever $\nabla(\ep\chi)\neq 0$) agrees with $J$ near $C_i$ for
$i>0$.  Further one easily sees that $\|J'-J\|_{C^1}=O(\ep)\leq
O(\delta^2)$.

\end{proof}

\begin{cor}\label{tvs}
Any set of distinct $J$-holomorphic curves $C_0,\ldots,C_m$ can be
perturbed to symplectic surfaces $C'_0,\ldots,C'_m$ whose
intersections are all transverse and positive, with $C'_i\cap
C'_j\cap C'_k=\varnothing$ when $i,j,k$ are all distinct.
Furthermore, there is an almost complex structure $J'$ arbitrarily
$C^1$-close to $J$ such that the $C'_i$ are $J'$-holomorphic.
\end{cor}
\begin{proof}
Assume that the process used in the proof of the above proposition
has been repeated to yield surfaces
$C^{\delta_0}_{0},\ldots,C^{\delta_i}_{i}$ each missing $p$ and
hitting the other $C_j$ transversely and positively.  Let our
neighborhood $U$ and the parameter $\delta_{i+1}$ be so small that
each $C^{\delta_{j}}_{j}$ ($j\leq i$) misses $B_{\delta_{i+1}}$
(this is possible since the $C^{\delta_{j}}_{j}$ all miss $p$); then
since $C_{i+1}^{\delta_{i+1}}\cap (X\setminus
B_{\delta_{i+1}})=C_{i+1}\cap (X\setminus B_{\delta_{i+1}})$, the
intersection points of $C^{\delta_{i+1}}_{i+1}$ with
$C_{j}^{\delta_j}$ ($j\leq i$) are the same as those of  $C_{i+1}$
and $C_{j}^{\delta_j}$, and so are transverse, positive, and away
from $p$.  By the proposition, we have the same conclusion for the
intersection points of $C^{\delta_{i+1}}_{i+1}$ with $C_j$
($j>i+1$).  So by induction we may perturb all of the $C_i$ to
$C'_i=C^{\delta_i}_{i}$ with the desired intersection configuration.
Moreover by choosing $\delta_0>\delta_1>\cdots \delta_{m-1}>0$ small
enough, the $C^{\delta_i}_{i}$ can be made arbitrarily $C^2$-close
to the $C_i$, so since the property of being a symplectic
submanifold persists under $C^1$-small perturbations, the
$C^{\delta_i}_{i}$ can be taken to all be symplectic.  Repeating
this local construction at all of the intersection points of two or
more of the $C_i$ gives the global result.
\end{proof}

\section{Towards a symplectic Nakai-Moishezon criterion}

In this subsection let $(M,J)$ be a minimal K\"ahler surface and
$H_J^{1,1}$ denote the real part of the $(1,1)-$subspace of
$H^2(M;{\mathbb C})$ determined by $J$. We apply Theorem \ref{main}
to study the symplectic classes in $H_J^{1,1}$.

Given a homology class $e$, we define the reflection along $e$ to be
$$R_e(\alpha)= \alpha-2{\alpha(e)\over e\cdot e}PD(e).$$
Notice that this is an automorphism of $H^2(M;{\mathbb Q})$
preserving the intersection form. But it is  an automorphism of
$H^2(M;{\mathbb Z})$ only if $e\cdot e=-1$ or $-2$. Geometrically,
the annihilator of $e$ is a hyperplane in $H^2(M;{\mathbb R})$ which
we call the ``$e$-wall,'' and $R_e$ is the reflection across this
hyperplane.

\begin{definition} A homology class $e$ is called small and effective if it is represented by a reduced irreducible holomorphic curve with
negative self-intersection.
\end{definition}

Notice that there is only one holomorphic curve $C$ representing a
small and effective class.

\begin{prop}\label{refl} Let $e$ be a small and effective class which is not represented by a curve of zero arithmetic genus and odd self-intersection. Then the  reflection of the K\"ahler chamber along the $e-$wall
is contained in the symplectic cone.
\end{prop}

\begin{proof} Let $x$ be a point in the K\"ahler cone. The K\"ahler cone is open in $H^{1,1}_{J}$, since
the sum of  a small closed real $(1,1)$ form and  a K\"ahler form on
a closed manifold is still a closed positive  $(1,1)$ form, hence a
K\"ahler form. Thus, for small $\epsilon$, $x-\epsilon e$ is also in
the K\"ahler cone, and hence represented by a K\"ahler form
$\omega$.  By Proposition \ref{pert}, we can perturb $C$ to get an
embedded $\omega-$symplectic surface, still denoted by $C$. Applying
Theorem \ref {main} to $\omega$ and $C$, we see that
$R_e(x)=[\omega_t]$ for some $t$.
\end{proof}

\begin{remark} For an embedded  $-2$ rational curve $C$, there is
 a diffeomorphism whose induced action on cohomology is $R_{[C]}$.
 Pulling the K\"ahler form back by this diffeomorphism gives an
 alternative way of enlarging the K\"ahler cone by reflection.
 However, this latter method, unlike Theorem \ref{main}, does not
 allow us to obtain symplectic forms in classes which vanish on
the $(-2)$-curve.
\end{remark}

We mention a simple case where the symplectic Nakai-Moishezon
criterion can be established.

\begin{prop} Suppose that $H_2(M;{\mathbb Z})$ contains only one small and effective class, $e$, and that $e$ is not represented by a sphere of odd square. Then every class $\alpha$ in the positive cone which is negative on $e$ lies in the image of the K\"ahler chamber under $R_e$.  Therefore the symplectic Nakai-Moishezon
criterion holds in this case.
\end{prop}
\begin{proof}  Suppose $e\cdot e=-k$ and $\alpha$ is as in the statement of the proposition.  Choose
$s>0$ such that
$$\alpha^2+2s|\alpha(e)|>s^2k>2s|\alpha(e)|.$$
Let $\beta=\alpha-sPD(e)$. Then
$$\beta(e)=\alpha(e)+sk>|\alpha(e)|, \beta^2=\alpha^2+2s|\alpha(e)|-s^2k>0, \beta\cdot \alpha=\alpha^2+|\alpha(e)|>0.$$
Therefore $\beta$ is in the K\"ahler cone by Theorem \ref{bl}. Now
apply Theorem \ref{main} to $\beta$.
\end{proof}

The much more common situation in which $M$ contains more than one
small and effective class is  more difficult to analyze.  We begin
by establishing the following finiteness result, which might be known to experts.

\begin{lemma} \label{finite} For any (1,1) class $\alpha$  with positive
square and in the positive cone, there are only finitely many
classes which are represented by reduced irreducible holomorphic
curves and pair non-positively with $\alpha$.  Further, the
intersection form on $M$ is negative definite on the subspace of
$H^{1,1}_{J}$ spanned by the Poincar\'e duals of these classes.
\end{lemma}
\begin{proof} Suppose $e_i$ are distinct such classes with negative
square which are represented by reduced irreducible holomorphic
curves. Notice that $e_i\cdot e_j\geq 0$ if $i\ne j$.

Then if a finite positive linear combination of $e_i$, say $\sum_i
a_ie_i$,   has non-negative square, it must be in the positive cone
or its boundary, as  $\omega$ is positive on each $e_i$, $a_i\geq
0$, and $\omega$ itself in the positive cone. By the Hodge index
theorem, as $\alpha$ is also in the positive cone, $\alpha$ is
strictly positive on $\sum_i a_ie_i$. But $\alpha$ is non-positive
on each $e_i$, so $\alpha$ is non-positive on $\sum_i a_ie_i$ as
$a_i\geq 0$.

This contradiction shows that any positive linear combination of the
$e_i$ has negative square. But this implies that for \emph{any}
$a_i\in {\mathbb R}$ not all zero we have, using positivity of
intersections between the distinct $e_i$,
\begin{eqnarray*} \left(\sum a_i e_i\right)^2 &=&\sum_i
a_{i}^{2}e_{i}^{2}+2\sum_{i<j} a_i a_je_i\cdot e_j\leq \sum_i
|a_{i}|^{2}e_{i}^{2}+2\sum_{i<j} |a_i| |a_j|e_i\cdot e_j\\
&\leq&\left(\sum |a_i| e_i\right)^2<0.\end{eqnarray*}  Thus the
$e_i$ are linearly independent, and they span a negative definite
subspace of $H_2(X;{\mathbb{Z}})$.
In particular, there are at most $h^{1,1}-1=b^-$ many $e_i$.
\end{proof}

In view of the lemma above, we make the following definition.  Here
${\mathcal P}$ denotes the positive cone in $H^{1,1}_{J}(X;{\mathbb
R})$.

\begin{definition} A finite set of small and effective classes $G=\{e_1,...,e_l\}$ is called
admissible if they are linearly independent, and the intersection
form on the subspace of $H_2(M;\mathbb{Z})$ generated by these $e_i$
is negative definite.
Given an admissible set $G$, the $G-$chamber
is
$${\mathcal C}(G)=\{\alpha\in{\mathcal P}|\alpha(e_i)\leq 0 \hbox{ if $e_i\in G$}, \quad  \alpha(e)> 0
\hbox{ if $e\not \in G$} \}.$$  The $G-$corner is
\[ {\mathcal C}^c(G)=\{\alpha\in{\mathcal P}|\alpha(e_i)= 0 \hbox{ if $e_i\in G$}, \quad  \alpha(e)> 0
\hbox{ if $e\not \in G$} \}.\]
\end{definition}

The following simple observation will be useful.

\begin{prop}\label{allplus} Let $M$ be a symmetric negative definite matrix such
that $M_{ij}\geq 0$ if $i\neq j$.  Then every entry of $-M^{-1}$ is
non-negative.\end{prop}
\begin{proof} By multiplying $M$ by a scalar assume without loss of generality that all diagonal
entries and all eigenvalues of $M$ are greater than $-1$.  Then,
where $I$ is the identity, $I+M$ has all its entries nonnegative and
all its eigenvalues between 0 and 1.  The latter condition implies
that we have a convergent Taylor series expansion \[
-M^{-1}=(I-(I+M))^{-1}=\sum_{n=0}^{\infty}(I+M)^n,\] and the
proposition follows from the fact that the set of matrices with all
entries nonnegative is closed under addition and multiplication.
\end{proof}

\begin{lemma} \label{chambers} The chambers ${\mathcal C}(G)$ for admissible sets $G$
form a partition of the positive cone and are all nonempty, as are
the $G-$corners $\mathcal{C}^c(G)$. Each $G-$chamber and each
$G-$corner is convex and hence connected.
\end{lemma}

\begin{proof} That the $\mathcal{C}(G)$ form a partition of the positive cone follows directly from Lemma \ref{finite}.
Convexity is obvious from the definitions.

To see that each $\mathcal{C}(G)\neq \varnothing$, let
$G=\{e_1,\ldots,e_n\}$ be an admissible set and denote by $M$ the
matrix representing the restriction of the intersection form to the
span of $G$, so that $M$ is negative definite.  Pick an arbitrary
$\alpha$ in the K\"ahler cone, and let $v_i=\langle
\alpha,e_i\rangle$, so that each $v_i>0$. Then where
$\vec{t}=-M^{-1}\vec{v}$ and $\alpha'=\alpha+\sum t_iPD(e_i)$, we
have $\langle \alpha',e_j\rangle=v_j-v_j=0$ for each $j$, and \[
(\alpha')^2=\alpha^2+2\vec{v}\cdot\vec{t}+(M\vec{t})\cdot\vec{t}=\alpha^2-(M\vec{t})\cdot\vec{t}\geq
\alpha^2>0\] since $M$ is negative definite, so $\alpha'$ is in the
positive cone.  Also, by Proposition \ref{allplus}, we have each
$t_i>0$ since each $v_i>0$, so if $e$ is small and effective with
$e\notin G$ then by positivity of intersections $\langle
\alpha',e\rangle\geq \langle \alpha,e\rangle>0$. Thus $\alpha'\in
\mathcal{C}^c(G)$, and $\mathcal{C}^c(G)$ is nonempty. Where
$s_i=-\sum (M^{-1})_{ij}$, $\alpha'+\ep\sum s_iPD(e_i)$ will
evaluate as $-\ep$ on each $e_i$, will be positive on each $e\notin
G$ (noting that each $s_i>0$), and will remain in the positive cone
for small $\ep>0$, so $\mathcal{C}(\{e_1,\ldots,e_n\})$ is also
nonempty.

\end{proof}

\begin{remark} By Theorem \ref{bl},
the K\"ahler cone is just ${\mathcal C}(\emptyset)$.  Within the
positive cone, the \emph{boundary} of the K\"ahler cone is the
disjoint union of the $\mathcal{C}^c(G)$ over the admissible sets
$G$.
\end{remark}

Applying Theorem \ref{main} with $\omega$ equal to a K\"ahler form,
$e=[C]$, and $t$ between $a/k$ and $2a/h$ shows that each chamber
${\mathcal C}(e)$ contains symplectic classes. Iterating Theorem
\ref{main}, the same can be said for any $G-$chamber ${\mathcal
C}(e_1,\ldots,e_n)$ with $e_i\cdot e_j=0$ for $i\ne j$.

We can apply Theorem \ref{configuration} to show that more general
$G-$chambers contain symplectic classes.  To do this, it suffices to
show that the corresponding $G$-corner contains symplectic classes,
since as in the proof of Lemma \ref{chambers} suitably chosen
arbitrarily small perturbations of these will lie in
$\mathcal{C}(G)$ and will remain symplectic.  Under suitable
hypotheses on the set $G$, we shall see that every class in the
$G$-corner $\mathcal{C}^c(G)$ contains symplectic forms.

Accordingly, let $\alpha\in H^{1,1}_{J}(M;\mathbb{R})$ be an
arbitrary class in the boundary of the K\"ahler cone and have
positive square, so that $\alpha$ satisfies $\langle
\alpha,D\rangle\geq 0$ for every effective divisor $D$. $\alpha$ is
then in some $G-$corner; say $G=\{e_1,\ldots,e_n\}$, so that
$\alpha$ vanishes only on the $e_i$ and the $PD(e_i)$ span a
negative definite subspace of $H^{1,1}_{J}(M;\mathbb{R})$.  Our
strategy for attempting to show that $\alpha$ contains symplectic
forms consists of two steps:
\begin{itemize}\item[(i)] Find $t_i>0$ such that $\alpha-\sum t_i PD(e_i)$ lies
in the K\"ahler cone. \item[(ii)]  Beginning with a K\"ahler form in
the class $\alpha-\sum t_i PD(e_i)$, apply the inflation procedure
sequentially to the $e_i$ (and/or smoothings of unions thereof) to
obtain a symplectic form in class $\alpha$. \end{itemize}

We shall show presently that step (i) can always be completed.
\begin{lemma}\label{bdrytoint}  If $\alpha$ and $e_i$ are as above, and if $s_i>0$ are
such that $\sum_i s_i e_i\cdot e_j<0$ for every $j$, then for $r>0$
sufficiently small, $\alpha-\sum rs_iPD(e_i)$ admits K\"ahler forms.
\end{lemma}
\begin{proof} Multiplying the $s_i$ by a small constant if necessary,
assume that $\beta:=\alpha-\sum s_i PD(e_i)$ is in the positive
cone. By Lemma \ref{finite}, there are then just finitely many
curves on which $\beta$ is non-positive; denote them by
$f_1,\ldots,f_m$ (note that the assumption on the $s_i$ implies that
none of the $f_j$ is among the $e_i$).  Now for each $f_j$ we have
$\langle \alpha,f_j\rangle>0$, so since there are only finitely many
$f_j$, for $r>0$ small enough $\alpha-\sum rs_i e_i=(1-r)c+rd$ will
be positive on each $f_j$. Meanwhile $\langle \alpha,e_i\rangle =0$
and $\langle \beta,e_i\rangle >0$, and if $C$ is any curve not among
the $e_i$ and $f_j$ both $\alpha$ and $\beta$ are positive on $[C]$,
so for $r>0$ $(1-r)\alpha+r\beta$ is also positive on all curves
other than those represented by the $f_j$. Hence by Theorem \ref{bl}
$(1-r)\alpha+r\beta$ admits K\"ahler forms if $r>0$ is small enough.
\end{proof}

\begin{cor} If $\alpha\in H^{1,1}_{J}$ has positive square and lies in
the boundary of the K\"ahler cone, and if $e_1,\ldots,e_n$ are the
homology classes of the finitely many curves on which $\alpha$
vanishes, then there are $t_i>0$ such that $\alpha-\sum t_i PD(e_i)$
contains K\"ahler forms.
\end{cor}
\begin{proof} By Lemma \ref{bdrytoint} it suffices to find $s_i>0$
such that $\sum_i s_i e_i\cdot e_j<0$ for every $j$; we then set
$t_i=rs_i$ for $r$ small.  Define the $n\times n$ matrix $M$ by
$M_{ij}=e_i\cdot e_j$.  $M$ is negative definite by Lemma
\ref{finite}, and its off-diagonal entries are nonnegative by
positivity of intersections, so $-M^{-1}$ has all nonnegative
entries by Proposition \ref{allplus}. Then for any $v_i>0$
($i=1,\ldots,n$), the $s_i=\sum -M^{-1}_{ik}v_k$ will each be
positive, and we have $\sum_i s_i e_i\cdot
e_j=-\sum_{i,k}M_{ji}M^{-1}_{ik}v_k=-v_j<0$, as desired.
\end{proof}

Carrying out step (ii) of our strategy is more difficult (and often
impossible).  As we allude to above, instead of applying inflation
sequentially to curves $C_i$ representing the $e_i$, we will
sometimes wish to smooth the union of the $C_i$ into an embedded
symplectic submanifold $C$ (as is always possible since the $C_i$
may be assumed to meet positively and transversely by Corollary
\ref{tvs}) and then apply the inflation procedure to $C$. Now $C$
will no longer be symplectic after we do this, and in the smoothing
construction $C$ will contain all but a small subset of each $C_i$,
so the $C_i$ won't be symplectic either.  As such, it will not be
possible to apply inflation to $C_i$ after we apply inflation to
$C$. The following trick allows us to evade this issue in certain
circumstances.

\begin{prop} \label{disjoin} Let $C_0,\ldots,C_k$ be symplectic surfaces such that
$C_0$ has only positive transverse intersections with the $C_i$
($i>0$).  Assume that \[ \#\left(C_0\cap\left(\bigcup_{i\geq 1}
C_i\right)\right)\geq -[C_0]^2\]  Then there exist symplectic
surfaces $\tilde{C_0}$ and $C$, homologous to $C_0$ and $\cup_{r\geq
0} C_r$ respectively, such that all intersections between
$\tilde{C_0}$ and $C$ are positive and transverse.\end{prop}

\begin{proof}
Where $m=-[C_0]^2$, assume that, for some points $p_1,\ldots,p_m$,
$C_0$ meets the surface $C_{i_j}$ at $p_j$; in complex coordinates
$(z,w)$ in a neighborhood $U_j$ around $p_j$ we may assume $C_0\cap
U_j=\{z=0\}$ and $C_{i_j}\cap U_j=\{w=0\}$. By exponentiating a
small scalar multiple of a smooth section of the normal bundle to
$C_0$ which vanishes negatively precisely at the $m=-[C_0]^2$ points
$p_j$, we take for $\tilde{C_0}$ a surface such that
$\tilde{C_0}\cap C_0=\{p_1,\ldots,p_m\}$ and, for each of the above
neighborhoods $U_j$, $\tilde{C_0}\cap U_j=\{(z,\epsilon \bar{z})\}$.
For $\epsilon$ small enough, $\tilde{C_0}$ will be sufficiently
$C^1$-close to $C_0$ as to guarantee that $\tilde{C_0}$ is
symplectic and (like $C_0$) only meets the $C_i$ ($i>0$) positively
and transversely. For $C$, we take a surface which coincides with
$\cup_{r\geq 0} C_r$ outside the $U_j$ and whose intersection with
$U_j$ is given by \[ C\cap U_j=\{(z,w)|zw=\delta f_{j}(z,w)\} \]
where $f_j$ is a \emph{real-valued} function supported on $U_j$ with
$f(p_j)\neq 0$ and $\delta$ is a complex constant chosen small
enough as to guarantee that $C$ is symplectic.  Now for any
$(z,w)\in \tilde{C_0}\cap U_j$, we have $zw\in \mathbb{R}\ep$, while
for any $(z,w)\in C\cap U_j$ we have $zw\in \mathbb{R}\delta$, so as
long as we choose $\epsilon,\delta\in \mathbb{C}$ to have different
phases we ensure that $C$ and $\tilde{C_0}$ have no intersections
within $\cup_{j\geq 1} U_j$.  By construction, any intersections of
$\tilde{C_0}$ with $C$ outside $\cup_{j\geq 1} U_j$ are positive and
transverse, proving the result.\end{proof}

There are many examples of intersection patterns of curves
$C_1,\ldots,C_n$ for which our methods give rise to symplectic
classes on the K\"ahler cone, but it does not seem possible at this
juncture to give a concise yet anywhere-near-exhaustive list of the
assumptions on the $[C_i]$ which are sufficient.  Instead, we shall
demonstrate the techniques on a particular complex surface, which we
believe illustrates nicely both the subtleties involved and the fact
that our construction gives rise to symplectic forms that cannot be
obtained by classical methods.

\subsection{The Kharlamov--Kulikov surface}

 If $(M,J)$ is a complex surface
admitting K\"ahler structures and $\mathcal{C}_J\subset H^{1,1}_{J}$
is the K\"ahler cone as given by the Buchdahl-Lamari theorem, then
every class in $\mathcal{C}_J+Re\,H^{2,0}_{J}$ is of course
represented by symplectic forms.  Although our method gives
seemingly new symplectic forms in classes $c$ outside
$\mathcal{C}_J$ in the presence of ($J$-holomorphic) curves of
negative square, a skeptic might imagine that if we were to vary the
complex structure on $M$ to some other (integrable) $J'$, then the
negative-square curves might disappear, and so these classes $c$
might lie in $\mathcal{C}_{J'}+Re\,H^{2,0}_{J'}$, in which case our
method would not have been necessary to obtain the new forms.

Now the list of underlying manifolds $M$ of complex surfaces for
which the effective cone is known for \emph{every} complex structure
on $M$ is rather short, so for most complex surfaces it is difficult
to tell whether our new forms could have been obtained by
algebro-geometric considerations.  In the case that $M$ is
\emph{rigid}, though, there is no room to vary $J$, and so we can
confidently assert that our main theorems give genuinely new forms
as soon as we know that there are curves of negative square in the
surface.  We present here an example of a rigid surface $K$,
borrowed from \cite{KK}, which contains several (21)  curves of
negative square intersecting each other in a nontrivial fashion, and
on which we can find symplectic forms in all classes in the positive
cone which are nonnegative on each of these 21 curves. It seems
likely (though we shall not attempt to prove this) that all curves
of negative square in $K$ lie in the cone generated by these 21
special curves; if this is indeed the case then it would follow that
the entire boundary of the K\"ahler cone of $K$ is contained in the
symplectic cone.  In any event, our results show that at least a
rather substantial portion of the boundary of the K\"ahler cone of
$K$ is contained in the symplectic cone, even though the standard
methods of K\"ahler geometry alone seem to provide no reason to
expect this to be the case.

We now recall the construction of $K$ from Section 2 of \cite{KK}.
Begin with an arbitrary smooth cubic curve in $\mathbb{C}P^2$, and
consider its 9 inflection points.  Since these inflection points are
each 3-torsion under the group law of the cubic, any line through
two of them also passes through a third which is distinct from the
first two; as such we obtain 12 lines each passing through precisely
3 of the inflection points. The dual arrangement provides us with 9
lines $L_1,\ldots,L_9$ and 12 points $p_{\{i,j,k\}}$ ($\{i,j,k\}\in$
 $\{\{1,2,3\},$ $\{1,4,7\},$ $\{1,5,9\},$ $\{1,6,8\},$ $\{2,4,9\},$ $\{2,5,8\},$ $\{2,6,7\},$ $\{3,4,8\},$
$\{3,5,7\},$ $\{3,6,9\},$ $\{4,5,6\},$ $\{7,8,9\}\}$) in (the dual
plane) $\mathbb{C}P^{2}$, with $p_{\{i,j,k\}}\in L_l$ iff
$l\in\{i,j,k\}$. Let $\sigma\co\tilde{\mathbb{P}}^2\to\mathbb{C}P^2$
denote the blowup at the various $p_{\{i,j,k\}}$; let
$E_{\{i,j,k\}}$ denote the corresponding exceptional divisors, and
let $L'_i$ denote the strict transform of $L_i$.  As is seen in
\cite{KK}, for suitable choices of a homomorphism $\phi\co
H_1(\tilde{\mathbb{P}}^2\setminus
\sigma^{-1}(\cup_{i=1}^{9}L_i);\mathbb{Z})\to
(\mathbb{Z}/5\mathbb{Z})^2$, the total space of the Galois cover
branched over $\cup_{i=1}^{9}L_i$ associated to $\phi$ will be
smooth.  Call this total space $K$ and the covering map $g\co K\to
\tilde{\mathbb{P}}^2$.

Write $C_i=g^{-1}(L'_i)$, $D_{\{i,j,k\}}=g^{-1}(E_{\{i,j,k\}})$.
Lemma 2.1 of \cite{KK} shows that each $C_i$ is a square-$(-3)$
curve of genus 4 and each $D_{\{i,j,k\}}$ is a square-$(-1)$ curve
of genus 2.  Further the canonical class of $K$ is ample and is
given by
\[ K_K=\frac{1}{3}PD(7\sum [C_i]+12\sum [D_{\{i,j,k\}}]);\] we have
$K_{K}^{2}=333$ and $e(K)=111$, so $K$ is the quotient of the unit
ball in $\mathbb{C}^2$ by a famous result of Miyaoka \cite{M} and
Yau \cite{Y}; a theorem of Siu \cite{Siu} then shows that $K$ is
rigid as promised.

\begin{theorem} Let $\alpha$ be any class in the positive cone of $H^{1,1}$
which is nonnegative on all holomorphic curves in $K$, and positive
on all curves whose homology classes are not in the cone spanned by
the $[C_i]$ and $[D_{\{i,j,k\}}]$. Then $\alpha$ is represented by
symplectic forms.\end{theorem}
\begin{proof}
First, note that the intersections of the distinct $C_i$ and
$D_{\{i,j,k\}}$ are given by
\[ [C_i]\cdot [C_j]=0;\quad [D_{\{i,j,k\}}]\cdot [D_{\{l,m,n\}}]=0; \quad
[C_{l}]\cdot [D_{\{i,j,k\}}]=\left\{\begin{array}{ll}1 & l\in\{i,j,k\}\\
0 & l\notin\{i,j,k\}\end{array}.\right.\]

Let $\Gamma$ denote the dual graph to the subset of
$\{[C_i],[D_{\{i,j,k\}}]\}$ on which $\alpha$ vanishes (in other
words, $\Gamma$ has a vertex for each element of this set, and the
number of edges connecting two distinct vertices of $\Gamma$ is the
intersection number of the corresponding pair of classes). If
$\Gamma$ were to contain a loop, then by virtue of the intersection
pattern of the $C_i$ and $D_{\{i,j,k\}}$ that loop would consist of
some number (say $a$) of curves
$A_0=C_{i_0},\ldots,A_{a-1}=C_{i_{a-1}}$ and an equal number of
curves $B_0=D_{\{i_0,j_0,k_0\}},\ldots,
B_{a-1}=D_{\{i_{a-1},j_{a-1},k_{a-1}\}}$ such that $[A_m]\cdot
[B_m]=[A_m]\cdot [B_{m+1}]=1$ for each $m$ (where $m\in
\mathbb{Z}/a\mathbb{Z}$). Hence since $[A_{m}]^{2}=-3$ and
$[B_{m}]^{2}=-1$,
\[ \left(\sum_{m=0}^{a-1}[A_m]+\sum_{m=0}^{a-1}[B_m]\right)^2\geq
-3a-a+2(2a)=0,\] which is impossible since $\alpha$ lies in the
positive cone and vanishes on $\sum [A_m]+\sum [B_m]$.

In general if $\Gamma$ contains a connected component with at least
3 distinct $[B_{m}]=[D_{\{i_m,j_m,k_m\}}]$ ($1\leq m\leq 3$), then
it contains a subgraph consisting of vertices
$\{[B_1],[C_{i_1}],[B_2],[C_{i_2}],[B_3]\}$ where $[B_1]\cdot
[C_{i_1}]=[B_2]\cdot [C_{i_1}]=[B_2]\cdot [C_{i_2}]=[B_3]\cdot
[C_{i_2}]=1$. But then
\[([B_1]+[C_{i_1}]+2[B_2]+[C_{i_2}]+[B_3])^2=-1-3-4-3-1+2+4+4+2=0,\] which is
again a contradiction since $\alpha$ is in the positive cone.
Likewise, if $\Gamma$ contains a connected component with three
distinct $[C_i]$ (say $[C_i]$, $[C_j]$, $[C_k]$), then it must also
contain some $[D_{\{i,j,l\}}]$ and $[D_{\{j,k,m\}}]$ and we see \[
([C_i]+3[D_{\{i,j,l\}}]+2[C_j]+3[D_{\{j,k,m\}}]+[C_k])^2=-3-9-12-9-3+6+12+12+6=0,\]
again a contradiction.

Now it will suffice to consider the case in which $\Gamma$ is
connected, since if it is not we can apply our argument successively
to each component.  Assuming $\Gamma$ is connected, then, the above
shows that it contains at most two $[C_i]$ and at most two
$[D_{\{i,j,k\}}]$, so after relabeling it is a subgraph of the graph
$\Gamma_0$ with vertices $[C_1]$, $[B_1]:=[D_{\{1,2,3\}}]$, $[C_2]$,
and $[B_2]:=[D_{\{2,4,9\}}]$, with just one edge each connecting
$[C_1]$ to $[B_1]$, $[B_1]$ to $[C_2]$, and $[C_2]$ to $[B_2]$.
Suppose that $\Gamma=\Gamma_0$.  Since $\alpha$ is positive on all
curves represented by classes which are not in the span of $[C_1]$,
$[B_1]$, $[C_2]$, and $[B_2]$, by taking $t>0$ small enough we
ensure that
\[ \alpha_0=\alpha-tPD(8[C_1]+21[B_1]+12[C_2]+14[B_2])\] will have the same property; we
calculate \[ \langle \alpha_0,[C_1]\rangle=3t,\,\langle
\alpha_0,[B_1]\rangle=t,\, \langle \alpha_0,[C_2]\rangle=t,\mbox{
and }\langle \alpha_0,[B_2]\rangle=2t,\] so $\alpha_0$ is
represented by K\"ahler forms. Apply Proposition \ref{disjoin}
twice: first to get a symplectic surface $\tilde{C}$ representing
$[C_1]+[B_1]$ and disjoint from a symplectic representative of
$[B_1]$, and then to get a symplectic surface $S$ representing
$[C_2]+[\tilde{C}]+[B_1]+[B_2]=[C_1]+2[B_1]+[C_2]+[B_2]$ which is
disjoint from $C_2$, $B_1$, and $B_2$.  $S$ then has positive genus
and square $-1$, so we can apply inflation to $S$ to get a
symplectic form in the class $\alpha_0+sPD[S]$ for any parameter $s$
less than $2\langle \alpha_0,[S]\rangle=16t$.  Take $s=8t$ to get a
symplectic form $\omega_1$ representing \[
\alpha_1=\alpha-tPD(5[B_1]+4[C_2]+6[B_2])\] with respect to which
$[B_1]$, $[C_2]$, and $[B_2]$ are symplectic. Now use Proposition
\ref{disjoin} to obtain a positive-genus $\omega_1$-symplectic
surface $S'$ representing $[B_1]+[C_2]+[B_2]$ and meeting $[B_1]$
and $[B_2]$ transversely and positively.  $[S']^2=-1$, and $\langle
\alpha_1,[S']\rangle=4t$, so inflation using $S'$ gives a symplectic
form $\omega_2$ in the class
\[ \alpha_1+4tPD[S]=\alpha-tPD([B_1]+2[B_2]).\]  Since $B_1\cdot B_2=0$, we can
now apply Theorem 1.1 rather directly to get the desired symplectic
form in $\alpha$, by first inflating using (say) $B_1$ and then
inflating using $B_2$.

In each case that $\Gamma$ is a \emph{proper} subgraph of
$\Gamma_0$, the desired symplectic representative of $\alpha$ can be
obtained by similar (but easier) arguments, which we leave to the
reader.\end{proof}

\subsection{A more general criterion}

As a more general example of the circumstances in which our methods
can be used to show that a class in the boundary of the K\"ahler
cone admits symplectic representatives, we present the following
theorem. Note that while condition (b) below is rather subtle,
condition (a) is occasionally easy to check; for instance it holds
for the canonical class in a minimal surface of general type and for
any class in the positive cone of a minimal surface of Kodaira
dimension 0 (though in both of these cases there exist other methods
to prove that such a class is in the symplectic cone).

\begin{theorem}\label{ADE}  Let $(M,\omega,J)$ be a K\"ahler surface and $\alpha\in
H^{1,1}_{J}$ any class in the positive cone such that
\begin{itemize} \item[(a)]  If $e\in H_2(M;\mathbb{Z})$ is
represented by a reduced, irreducible holomorphic curve of negative
square, then $\langle \alpha,e\rangle\geq 0$, with equality only if
$e^2=-2$ or $e^2=-1$ and $g(e)>0$; and \item[(b)]  There are no
$E_6$-trees of holomorphic curves of square $-2$ on which $a$
vanishes.\end{itemize}  Then $\alpha$ is represented by symplectic
forms deformation equivalent to $\omega$.
\end{theorem}
\begin{proof} (Sketch) Using negative-definiteness as in the case of
the Kharlamov--Kulikov surface, one first shows that each connected
component of the dual graph of the curves on which $\alpha$ is
negative either \begin{itemize}\item contains just one curve of
square $-1$ and (say) $n-1$ curves of square $-2$, in which case the
dual graph is the Dynkin diagram $A_n$, with the square-$(-1)$ curve
as one of the univalent vertices; or
\item consists entirely of square-$(-2)$ curves, in which case it is
one of the ADE Dynkin diagrams.\end{itemize} Now assumption (b) in
the statement of the theorem restricts the Dynkin diagrams that can
appear to $A_n$ and $D_n$, and is imposed because our methods do not
seem strong enough to apply to the cases of $E_6$, $E_7$, or $E_8$.
In the cases of $A_n$ and $D_n$, an approach parallel to that used
in the case of the Kharlamov--Kulikov surface provides the desired
form; the details of this are left as a mildly amusing exercise to
the interested reader.\end{proof}

\end{document}